\documentclass[twoside,leqno,10pt]{amsart}
\usepackage{amsfonts}
\usepackage{amsmath}
\usepackage{amscd}
\usepackage{amssymb}
\usepackage{amsthm}
\usepackage{latexsym}
\usepackage{bbm}
\numberwithin{equation}{section}

\begin{document}
\newtheorem{theorem}{Theorem}
\newtheorem{lemma}{Lemma}
\newtheorem{prop}{Proposition}
\newtheorem{corollary}{Corollary}
\newtheorem{conjecture}{Conjecture}
\numberwithin{equation}{section}
\newcommand{\dif}{\mathrm{d}}
\newcommand{\intz}{\mathbb{Z}}
\newcommand{\ratq}{\mathbb{Q}}
\newcommand{\natn}{\mathbb{N}}
\newcommand{\comc}{\mathbb{C}}
\newcommand{\rear}{\mathbb{R}} 
\newcommand{\prip}{\mathbb{P}}
\newcommand{\uph}{\mathbb{H}}
\newcommand{\fie}{\mathbb{F}}

\title{The Sato-Tate Conjecture on Average for Small Angles}
\date{\today}
\author{Stephan Baier \and Liangyi Zhao}
\maketitle

\begin{abstract} 
We obtain average results on the Sato-Tate conjecture for elliptic curves for small angles. 
\end{abstract}

\noindent {\bf Mathematics Subject Classification (2000)}: 
11G05\newline

\noindent {\bf Keywords}: Sato-Tate conjecture, average Frobenius distribution

\section{Introduction and main results}
Let $E$ be an elliptic curve over $\ratq$. For any prime number $p$ of good reduction, let $\lambda_E(p)$ be the trace of the Frobenius morphism of $E/\fie_p$. Then the number of points on the reduced curve modulo $p$ equals $\#E(\fie_p)= p+1- \lambda_E(p)$. By Hasse's theorem, there exists a unique angle $0\le \theta\le \pi$ such that
$$\lambda_E(p)=\sqrt{p}\left(e^{i\theta_E(p)}+e^{-i\theta_E(p)}\right)=
2\sqrt{p}\cos \theta_E(p).$$ 
It is natural to ask how $\theta_E(p)$ varies with $p$. \newline

If $E$ has complex multiplication(CM), the answer is easy. By Deuring's work \cite{Deur} half of the primes satisfy $\lambda_E(p)=0$ (these primes are called ``supersingular primes''), and for the remaining half of the primes the angles $\theta_E(p)$ are distributed uniformly in the interval $[0,\pi]$.  The reason is that, in the CM-case, the primes $p$ with $\theta_E(p)$ in a fixed range are given by $p=f(u,v)$, where $f(u,v)$ is a certain positive definite binary quadratic form, with $|u|/p$ in the corresponding range. The uniform distribution follows now from the work of Hecke \cite{Heck1}, \cite{Heck2} on prime ideals of imaginary quadratic number fields in sectors. \newline

The situation is more complicated when $E$ does not have complex multiplication. For this case, Sato and Tate \cite{Tate} formulated the following conjecture.\\ \\
{\bf Sato-Tate Conjecture:} Suppose $E$ is an elliptic curve over $\ratq$ which
does not admit complex multiplication. For any $0\le \theta_1\le \theta_2\le
\pi$, and $x\ge1$, let
$$
\pi_E^{\theta_1,\theta_2}(x):=|\{p\le x\ :\ \theta_1\le \theta_E(p)\le
\theta_2\}|.
$$
Then 
$$
\lim\limits_{x\rightarrow\infty} \frac{\pi_E^{\theta_1,\theta_2}(x)}{\pi(x)}
=\frac{2}{\pi}\int\limits_{\theta_1}^{\theta_2} \sin^2\theta\ {\rm d}\theta,
$$
where $\pi(x)$ is the number of primes not exceeding $x$. \newline

Let 
$$
\Theta_E^{\theta_1,\theta_2}(x):=\sum\limits_{\substack{p\le x\\ 
\theta_1\le \theta_E(p)\le \theta_2}} \log p.
$$
Then the Sato-Tate conjecture is equivalent to 
\begin{equation} \label{ST}
\lim\limits_{x\rightarrow\infty} \frac{\Theta_E^{\theta_1,\theta_2}(x)}{x}
=\frac{2}{\pi}\int\limits_{\theta_1}^{\theta_2} \sin^2\theta\ {\rm d}\theta.
\end{equation}
By a change of variables in \eqref{ST}, the conjecture can be stated in
the following equivalent way. \\ \\
{\bf Sato-Tate Conjecture:} 
Suppose $E$ is an elliptic curve over $\ratq$ which does not admit complex multiplication. For any $-1\le \alpha\le \beta\le 1$, and $x\ge 1$, let
\[ \Theta_E(\alpha,\beta;x):=\sum\limits_{\substack{p\le x\\ \alpha\le \lambda_E(p)/(2\sqrt{p})\le\beta}}\log p. \]
Then 
\[ \lim\limits_{x\rightarrow\infty} \frac{\Theta_E(\alpha,\beta;x)}{x} =\frac{2}{\pi}\int\limits_{\alpha}^{\beta} \sqrt{1-t^2}\ {\rm d}t. \]

In \cite{CHT}, \cite{HSBT} and \cite{Tayl}, L. Clozel, M. Harris, N. Shepherd-Barron and R. Taylor have proved the Sato-Tate conjecture for all elliptic curves $E$ over totally real fields (in particular, over $\mathbbm{Q}$) satisfying the mild condition of having multiplicative reduction at some prime. In the following, we say some few words about the background of their proof.

It can be shown that the Sato-Tate conjecture holds true if and only if the matrices
$$
\left( \begin{array}{*{2}{c}} e^{i\theta_E(p)} & 0 \\ 0 & e^{-i\theta_E(p)}\\ \end{array} \right)
$$
are uniformly distributed with respect to the Haar measure in the space
of conjugacy classes of the compact group SU$_2(\comc)$ of unitary
$2\times 2$ matrices over $\comc$ with determinant 1.  Building on this observation, Langlands \cite{Lang} and Serre in the letter at the end
of \cite{Sh} related the Sato-Tate conjecture to
symmetric power $L$-functions associated to the curve $E$. These $L$-functions are (in normalized form) defined by
$$
L_m(E;s):=\prod\limits_{p\nmid \Delta(E)}\ \prod\limits_{j=0}^m
\left(1-e^{(m-2j)i\theta_E(p)} p^{-s}\right)^{-1},
$$
for $m=1,2,...$,  where $\Delta(E)$ is the discriminant of
$E$. $L_m(E;s)$ converges absolutely for $\Re s>1$. Langlands proved
that the Sato-Tate conjecture holds if for all non-negative integers
$m$, $L_m(E;s)$ can be extended to an analytic function on $\Re s\ge 1$
and does not vanish on the line $\Re s=1$. By an argument given in
\cite{MuMu}, the existence of an analytic continuation of $L_m(E;s)$ to
the closed half plane $\Re s\ge 1$ implies the non-vanishing of
$L_m(E;s)$ on the line $\Re s=1$.  Langlands conjectured that the above symmetric power $L$-functions
extend to an entire function and coincide with certain automorphic
$L$-functions.  This conjecture implies the Sato-Tate conjecture by the
above consideration.  See \cite{Murt} for an overview of these results. \newline

Recently, the authors of \cite{CHT}, \cite{HSBT} and \cite{Tayl} have proved that, under the above-mentioned mild restriction on $E$, these symmetric power
$L$-functions have meromorphic continuation to the whole complex plane,
satisfy the expected functional equation and are holomorphic and non-zero
for $\Re s \geq 1$.  These conditions suffice to prove the Sato-Tate
conjecture.  We note here that our results, which we shall presently
state, do not follow from those in \cite{CHT}, \cite{HSBT} and \cite{Tayl} due to the lack of information about the order of magnitude of the error term in the Sato-Tate conjecture. In particular, as far as we know, there are no results in the literature about uniformity of this error term
with respect to the coefficients $a$ and $b$ of the elliptic curves. \newline

To establish an asymptotic estimate for $\pi_E^{\theta_1,\theta_2}(x)$ or $\Theta_E^{\theta_1,\theta_2}(x)$ if the angle $[\theta_1,\theta_2]$ is very small (that is, $\theta_2-\theta_1<x^{-\delta}$ for some positive $\delta$), one would need zero density estimates for these symmetric power $L$-functions which do not depend too much on $m$.  It seems that there are no zero density estimates in the literature which are sufficient for this purpose. In this paper we prove that the expected asymptotic estimate for small angles 
holds on average over elliptic curves 
$$
E(a,b) \ : \ Y^2=X^3+aX+b
$$
with $\vert a\vert \le A$, $\vert b\vert \le B$. For the sake of 
technical simplicity, we work with $\Theta_{E(a,b)}(\alpha,\beta;x)$ 
with a 
small interval $[\alpha,\beta]$ instead of 
$\Theta_{E(a,b)}^{\theta_1,\theta_2}(x)$
with a small angle $[\theta_1,\theta_2]$. We further suppose without 
loss of generality that $\alpha$ and $\beta$ are positive. The case when
$\alpha$ and $\beta$ are negative can be treated in a similar way. The 
distribution of supersingular primes $p$ with $\lambda_p=0$ was investigated in
\cite{FoMu}.\newline

In the sequel, we set
\begin{equation}\label{Fdef}
F(\alpha,\beta):=\frac{2}{\pi}
\int\limits_{\alpha}^{\beta} \sqrt{1-t^2}\ {\rm d}t \; \mbox{ and } \;
\gamma:=\beta-\alpha. 
\end{equation}
Moreover, following the general convention, we use $\varepsilon$ to denote a small positive constant which may not be the same at each occurrence.
Our first result is:

\begin{theorem} \label{TH1}
Let $x\ge 1$, $0< \alpha\le \beta\le 1$ and $A,B\ge 1$. 
Assume that $x^{\varepsilon-5/12}\le \gamma/\beta \le x^{-\varepsilon}$ and $F(\alpha,\beta)\ge x^{-1/2+\varepsilon}$. 
Then, for every $c>0$, we have
\[ \frac{1}{4AB}\sum\limits_{|a|\le A} \sum\limits_{|b|\le B} \Theta_{E(a,b)}(\alpha,\beta;x) = xF(\alpha,\beta) +O\left(\frac{xF(\alpha,\beta)}{\log^c x}+x^{3/2+\varepsilon}
\left(F(\alpha,\beta)
\left(\frac{1}{A}+\frac{1}{B}\right)+
\frac{F(\alpha,\beta)^{1/2}}{(AB)^{1/2}}\right)\right), \]
where the implied $O$-constant depends only on $c$.
\end{theorem}
 
From Theorem~\ref{TH1}, we immediately deduce the following average result for small intervals $[\alpha,\beta]$.

\begin{corollary} \label{CO1}
Let $x\ge 1$, $0< \alpha\le \beta\le 1$. Assume that 
$x^{\varepsilon-5/12}\le \gamma/\beta \le x^{-\varepsilon}$ and $F(\alpha,
\beta)\ge x^{-1/2+\varepsilon}$.  
Then, if $A,B>x^{1/2+\varepsilon}$ and 
$AB>x^{1+\varepsilon}/F(\alpha,\beta)$, we have, for every $c>0$,
\[ \frac{1}{4AB}\sum\limits_{|a|\le A} \sum\limits_{|b|\le B} \Theta_{E(a,b)}(\alpha,\beta;x) =xF(\alpha,\beta) \left(1+O\left(\frac{1}{\log^c x}\right)\right), \]
where the implied $O$-constant depends only on $c$.
\end{corollary}

We note that under the assumptions in Corollary~\ref{CO1} the contribution of CM-curves among the curves $E(a,b)$ is negligible since it is well-known that the number of these curves with $|a| \leq A$ and $|b| \leq B$ under consideration is $O(A+B)$. Therefore Theorem~\ref{TH1} and Corollary~\ref{CO1} are consistent with the Sato-Tate Conjecture. \newline
 
Moreover, we shall establish the following second moment estimate.

\begin{theorem} \label{TH2}
Let $x\ge 1$, $0< \alpha\le \beta\le 1$. 
Assume that $x^{\varepsilon-5/12}\le \gamma/\beta \le x^{-\varepsilon}$, 
$F(\alpha,\beta)\ge x^{-1/2+\varepsilon}$,  
$A,B>x^{1/2+\varepsilon}$ and 
$AB>x^{1+\varepsilon}/F(\alpha,\beta)$. Then, for every $c>0$, we have
\begin{equation} \label{T3}
\begin{split}
\frac{1}{4AB}& \sum\limits_{|a|\le A} \sum\limits_{|b|\le B} \left\vert \Theta_{E(a,b)}(\alpha,\beta;x)-xF(\alpha,\beta)\right\vert^2\\ 
\ll& \frac{x^2F(\alpha,\beta)^2}{\log^{c} x}+ x^{\varepsilon}  \Bigg(xF(\alpha,\beta)\log\log(10AB) + x^3F(\alpha,\beta)^2\left(\frac{1}{A}+\frac{1}{B}\right)+ \frac{x^3F(\alpha,\beta)}{(AB)^{1/2}}\Bigg),
\end{split}
\end{equation}
where the implied $\ll$-constant depends only on $c$.
\end{theorem}
 
Theorem~\ref{TH2} implies that the expected asymptotic estimate holds 
for almost all elliptic curves in a sufficiently
large box even if the interval $[\alpha,\beta]$ is very small.

\begin{corollary} \label{CO2}
Let $x\ge 1$, $0< \alpha\le \beta\le 1$. Assume that 
$x^{\varepsilon-5/12}\le \gamma/\beta \le x^{-\varepsilon}$, 
$F(\alpha,\beta)\ge x^{-1/2+\varepsilon}$, 
$A,B>x^{1+\varepsilon}$, $x^{2+\varepsilon}/F(\alpha,\beta)^2<AB<
\exp(\exp(x^{1-\varepsilon}))$. 
Then, for any $c,d>0$, we have
$$
\Theta_{E(a,b)}(\alpha,\beta;x)
=xF(\alpha,\beta)
\left(1+O\left(\frac{1}{\log^c x}\right)\right)
$$
for all $|a|\le A$, $|b|\le B$ with at most $O\left(AB/(\log x)^d\right)$ 
exceptions. Here the implied $O$-constant depends only on $c$ and $d$.
\end{corollary}

Our method is essentially a combination of those used  in \cite{Baie} and \cite{DaPa}.  However, we need to modify these methods appropriately to make them applicable to short intervals. Parts of our section 4 follow closely to the work of David-Pappalardi in \cite{DaPa} and parts of sections 5, 6 and 7 follow closely to the work of Baier in \cite{Baie}.  We included these computations for the sake of completeness. \newline

Corollaries~\ref{CO1} and~\ref{CO2} improve some results in a recent
preprint of James and Yu \cite{JaYu} on the Sato-Tate Conjecture on
average. They dealt only with fixed intervals $[\alpha,\beta]$, whereas
in our present paper the length of this interval may depend on the
parameter $x$. Moreover, our method is different from that in
\cite{JaYu}. James and Yu used the circle method which is avoided in the
present paper.  It is also note-worthy that Banks and Shparlinski
\cite{WDBIES} have recently improved our results in certain cases using very
different methods while in other cases our results are better.

\section{General approach}
As in \cite{DaPa} and \cite{Baie}, the following observations are our 
starting point.  

\begin{lemma} [Deuring] \label{LE1} 
For $0<r< 2\sqrt{p}$, the number of $\fie_p$-isomorphism classes of 
elliptic curves over $\fie_p$ with $p+1-r$ points is the total number of 
ideal classes of
the ring $\intz[(D+\sqrt{D})/2]$, where $D=r^2-4p$ is a negative integer which
is congruent to $0$ or $1$ modulo $4$. This number is the Kronecker class 
number $H(r^2-4p)$.
\end{lemma} 

\begin{lemma} \label{LE2} 
Suppose that $p\not=2,3$. Then any elliptic curve over $\fie_p$ has a model
$$
E\ :\ Y^2=X^3+aX+b
$$
with $a,b\in \fie_p$. The elliptic curves $E'(a',b')$ over $p$, which are
$\fie_p$-isomorphic to $E$, are given by all the choices
$$
a'=\mu^4a \ \ \ \ \mbox{\rm and } \ \ \ b'=\mu^6b
$$
with $\mu\in \fie_p^*$. The number of such $E'$ is
\begin{eqnarray*}
(p-1)/6, & &  \mbox{\rm if } a=0\ \mbox{\rm and } p\equiv 1\ \pmod{3};\\
(p-1)/4, & &  \mbox{\rm if } b=0\ \mbox{\rm and } p\equiv 1\ \pmod{4};\\
(p-1)/2, & &  \mbox{\rm otherwise.}
\end{eqnarray*}
\end{lemma}

Lemmas~\ref{LE1} and~\ref{LE2} imply that the number of curves $E(a,b)$ with $a,b\in \intz$, $0\le a,b<p$ and $\lambda_{E(a,b)}(p)=r$ is
\begin{equation} \label{NC}
\frac{p H(r^2-4p)}{2}+O(p).
\end{equation}
Now we write
\begin{equation} \label{goal}
\frac{1}{4AB} \sum\limits_{|a|\le A} \sum\limits_{|b|\le B} \Theta_{E(a,b)}(\alpha,\beta;x) = \frac{1}{4AB}\sum\limits_{p\le x}\log p \sum\limits_{2\sqrt{p}\alpha\le r\le 2\sqrt{p}\beta} \sharp\{|a|\le A,\ |b|\le B\ :\ \lambda_p(E(a,b))=r\}.
\end{equation}
Using \eqref{NC}, the term on the right-hand side of \eqref{goal} is
\begin{equation} \label{DaPa}
\frac{1}{4AB} \sum\limits_{p\le x}\log p \sum\limits_{2\sqrt{p}\alpha\le r\le 2\sqrt{p}\beta} \left(\frac{2A}{p}+O(1)\right) \left(\frac{2B}{p}+O(1)\right)\left(\frac{p H(r^2-4p)}{2}+O(p)\right).
\end{equation}
Here the error terms are estimated trivially.  A more precise evaluation
of the error term will be carried out in section 6.  The main term in \eqref{DaPa} is 
\begin{equation} \label{Main}
\sum\limits_{p\le x}\log p \sum\limits_{2\sqrt{p}\alpha\le 
r\le 2\sqrt{p}\beta} \frac{H(r^2-4p)}{2p}.
\end{equation}
We shall derive an asymptotic estimate for this main term in section 4. 

\section{Preparations for the estimation of the main term}
In this section we assemble some lemmas which we shall later need for the 
estimation of the main term. First, we express $H(r^2-4p)$ in terms 
of the value of a certain Dirichlet $L$-function at the point $s=1$. 
For any negative integer $d$ let $\chi_d$ be the character modulo $|d|$ given
by 
$$
\chi_d(n)=\left(\frac{d}{n}\right),
$$
where the right-hand side is the Jacobi symbol. Furthermore, let
$L(s,\chi)$ be the Dirichlet $L$-function associated to a Dirichlet character
$\chi$. Then we have the following.

\begin{lemma} \label{Hlemma} If $0<r< 2\sqrt{p}$, then
$$
H(r^2-4p)=\frac{1}{\pi}\sum\limits_{\substack{f,d\\ r^2-4p=df^2\\ d\equiv 0,1\
\bmod{4}}} \sqrt{|d|}L(1,\chi_d). 
$$
\end{lemma}
\begin{proof} This follows from (23) and (30) in \cite{DaPa}. \end{proof}

We shall also need the classical P\'olya-Vinogradov inequality for character sums. \newline

\begin{lemma} \label{PolyaVinolemma}
Let $q,N\in \mathbbm{N}$, $q\ge 2$ and $\chi$ be any non-principal character 
modulo $q$. Then
$$
\sum\limits_{n\le N} \chi(n)\ll \sqrt{q}\log q.
$$
\end{lemma}
\begin{proof} This is the well-known inequality of P\'olya-Vinogradov and can be found as Theorem 12.5 in \cite{IwKo}. \end{proof}

In the sequel, for $x,y\ge 1$ and $(a,q)=1$, we set
\begin{equation} \label{defpsi1}
\Theta(x,y;q,a):=\sum\limits_{\substack{x<p\le x+y\\ p \equiv a\
\bmod{q}}} \log p, \; \mbox{and} \; E(x,y;q,a):=\Theta(x,y;q,a)-\frac{y}{\varphi(q)}.
\end{equation}
We shall need the following estimate of $E(x,y;q,a)$ for small moduli $q$. 

\begin{lemma} \label{arithlemma}
Let $c_1,c_2,\varepsilon>0$ be given. 
Then, if $x\ge 2$, 
$x^{7/12+\varepsilon}\le y\le x$, $q\le (\log x)^{c_1}$ and 
$(a,q)=1$, we have
\begin{equation} \label{Estimate}
E(x,y;q,a)\ll \frac{y}{(\log x)^{c_2}},
\end{equation}
where the implied $\ll$-constant depends only on $c_1$, $c_2$ and 
$\varepsilon$.
\end{lemma}
\begin{proof} For $q=1$ this result is well-known. Its proof uses density estimates for the zeros of the Riemann zeta function as well as Vinogradov's zero-free region (see the proof of Theorem 12.8. in \cite{Ivic}). In a similar manner one can prove \eqref{Estimate} for any $q\le (\log x)^{c_1}$ by writing $E(x,y;q,a)$ as a linear combination of character sums and using density estimates for the zeros of the corresponding Dirichlet $L$-functions (see section 10.4 in \cite{IwKo}) as well as an analogue of Vinogradov's zero-free region for these $L$-functions which can be established by a similar method as for the Riemann zeta function (see section 8.5 in \cite{IwKo}). \end{proof}

We also need the following classical large sieve inequality for Dirichlet characters.

\begin{lemma} [Large Sieve] \label{lslemma}
Let $\{ a_n \}$ be a sequence of complex numbers.  Suppose that $M \in \intz$, $N, Q \in \natn$. Then we have
\begin{equation} \label{classlseq}
\sum_{q=1}^Q \frac{q}{\varphi(q)} \sideset{}{^{\star}}\sum_{\chi \bmod q} \left| \sum_{n=M+1}^{M+N} a_n \chi (n) \right|^2 \leq (Q^2 +N ) \sum_{n=M+1}^{M+N} |a_n|^2,
\end{equation}
where $\sideset{}{^{\star}}\sum$ henceforth denotes the sum over primitive characters to the specified modulus.
\end{lemma}
\begin{proof} See for example \cite{Daven} for the proof. \end{proof}

We are now ready to state the following generalization of the Barban-Davenport-Halberstam Theorem for short intervals. 

\begin{lemma} \label{BDVlemma}
Let $c_3,\varepsilon>0$ be given. Then, if $x\ge 2$, $x^{7/12+\varepsilon}\le y\le x$ and $1\le Q\le y$, we have
\begin{equation} \label{BDV}
\sum\limits_{q\le Q} \sum\limits_{\substack{a=1\\ (a,q)=1}}^q |E(x,y;q,a)|^2\ll Qy(\log x)^2 + \frac{y^2}{(\log x)^{c_3}},
\end{equation}
where the implied $\ll$-constant depends only on $c_3$ and $\varepsilon$.
\end{lemma}
\begin{proof} For $y=x$ this is the well-known classical Barban-Davenport-Halberstam Theorem (see section 17.4 in \cite{IwKo}). Its proof uses the Siegel-Walfisz Theorem for small moduli $q\le (\log x)^A$ and the large sieve, Lemma \ref{lslemma}, for the sum over the remaining moduli $q>(\log x)^A$.  To establish \eqref{BDV} for any $y$ in the range $x^{7/12+\varepsilon}\le y\le x$, one may use a similar method, where the Siegel-Walfisz Theorem is replaced by Lemma~\ref{arithlemma}.
\end{proof}

Furthermore, as in \cite{DaPa}, we set
\begin{equation} \label{constantfrn}
c_f^r(n):=\sum\limits_{\substack{a=1\\ (a,4n)=1\\ (r^2-af^2,4n)=4}}^{4n}
\left(\frac{a}{n}\right)
\end{equation}
and 
\begin{equation} \label{LTconstant}
K_r:=\prod\limits_{l|r} \left(1-\frac{1}{l^2}\right)^{-1}
\prod\limits_{l\nmid r} \frac{l(l^2-l-1)}{(l-1)(l^2-1)},
\end{equation}
where the first product on the right-hand  side runs over all prime divisors 
$l$ of $r$, and the second product runs over all primes $l$ which do not
divide $r$. We shall need the following asymptotic estimate for the estimation
of the main term in section 4. 

\begin{lemma} \label{maintermlemma}
Let $U,V\ge 1$ and $r>0$. Then
$$
\sum\limits_{n\le U}\sum\limits_{\substack{f\le V\\ (2r,f)=1}}
\frac{c_f^r(n)}{fn\varphi\left(nf^2\right)}=
K_r+O\left(\frac{1}{\sqrt{U}}+\frac{1}{V^2}\right),
$$
with an absolute $O$-constant.
\end{lemma}
\begin{proof} This follows from (22) and Lemma 4.1 in \cite{DaPa}. We note that these estimates do not depend on $r$.\end{proof}

Finally, we shall need the following estimate for $K_r$ on average. 

\begin{lemma} \label{LTconstantlemma}
Let  $u\ge 0$ and $v\ge 1$. Then, for $K_r$ being defined as in 
\eqref{LTconstant}, we have
\begin{equation} \label{Kaverage}
\sum\limits_{u<r\le u+v} K_r = v+O(1), 
\end{equation}
with an absolute $O$-constant.
\end{lemma}
\begin{proof} Let 
$$
C=\prod\limits_{l} \left(1+\frac{1}{l(l^2-l-1)}
\right)^{-1} \ \ \mbox{ and } \ \
f(r):=\prod\limits_{l|r} \left(1+\frac{1}{l^2-l-1}\right).
$$
Then 
\begin{equation} \label{Kidentity}
K_r=f(r)C.
\end{equation}
In the following, we estimate $f(r)$ on average over intervals. To this end,
we define 
the arithmetic function $g$ to be the Dirichlet convolution of $f$ and 
the M\"obius function $\mu$, that is, $g=f\ast \mu$. Then it is easy to
check that 
$$
g(n)=\frac{\mu(n)^2}{\prod\limits_{l|n} (l^2-l-1)}.
$$
Using this and the M\"obius inversion formula, we deduce 
\begin{equation} \label{faverage}
\sum\limits_{u<r\le u+v} f(r) = \sum\limits_{u<r\le u+v} \sum\limits_{d\vert r} g(d) = \sum\limits_{d=1}^{\infty} g(d)\left(\frac{v}{d}+O(1)\right) = v\sum\limits_{d=1}^{\infty} \frac{g(d)}{d} +O(1) = vC^{-1}+O(1).
\end{equation}
Now \eqref{Kaverage} follows from \eqref{Kidentity} and \eqref{faverage}. 
\end{proof}

\section{Estimation of the main term}
We are now ready to prove the following asymptotic estimate for the main term
in \eqref{Main}.
\begin{theorem} \label{maintermtheorem}
Let $x\ge 2$ and $0< \alpha\le \beta\le 1$. 
Suppose that $x^{2\varepsilon-5/12}
\le \gamma/\beta\le 
x^{-\varepsilon}$ and $\alpha\ge x^{-1/2+\varepsilon}$. Then, for any
$c>0$, we have
$$
\sum\limits_{p\le x}\log p \sum\limits_{2\sqrt{p}\alpha\le 
r\le 2\sqrt{p}\beta} \frac{H(r^2-4p)}{2p} = 
xF(\alpha,\beta)\left(1+O \left( \frac{1}{\log^c x}\right) \right).
$$  
\end{theorem}

\begin{proof} We have 
\begin{equation} \label{Homer}  
\begin{split}
\sum\limits_{p\le x} \log p \sum\limits_{2\sqrt{p}\alpha\le 
r\le 2\sqrt{p}\beta} \frac{H(r^2-4p)}{2p} = & \sum\limits_{1\le r\le x^{2\varepsilon}\alpha\beta^{6/5}/ \gamma^{6/5}} \ \sum\limits_{r^2/(2\beta)^2\le p\le r^2/(2\alpha)^2}\frac{H(r^2-4p)}{2p}\log p \\
& \hspace*{.25in} + \sum\limits_{x^{2\varepsilon}\alpha\beta^{6/5}/ \gamma^{6/5} < r \le 2\sqrt{x}\alpha} \ \sum\limits_{r^2/(2\beta)^2\le p\le r^2/(2\alpha)^2}\frac{H(r^2-4p)}{2p}\log p \\
& \hspace*{.25in} +\sum\limits_{ 2\sqrt{x}\alpha< r\le 2\sqrt{x}\beta} \ \sum\limits_{r^2/(2\beta)^2\le p\le x}\frac{H(r^2-4p)}{2p}\log p, \\
& = S_1+S_2+S_3, \; \mbox{say},
\end{split}
\end{equation}
where we have taken into account the condition $\gamma/\beta \geq x^{2\varepsilon - 5/12}$.We shall deal with $S_1$ and $S_3$ later.  Using Lemma \ref{Hlemma}, we may write $S_2$ in the form
\begin{equation} \label{Burns}
\frac{1}{2\pi} 
\sum\limits_{x^{2\varepsilon}\alpha\beta^{6/5}/ \gamma^{6/5} \leq r\le 2\sqrt{x}\alpha} \sum_{f \leq 2\sqrt{x}} \frac{1}{f} 
\sum\limits_{p \in 
S_{f,r}(r^2/(2\alpha)^2)} \frac{\sqrt{4p-r^2}}{p} L(1,\chi_d) \log p, 
\end{equation}
where 
\[ S_{f,r}(t) = \left\{ \left( \frac{r}{2 \beta} \right)^2 \leq p \leq t : f^2|(4p-r^2), d = \frac{r^2-4p}{f^2} \equiv 0,1 \bmod{4} \right\}. \]
We now develop an asymptotic estimate for the inner-sums over $f$ and $p$ of \eqref{Burns} for every $r$.  It suffices to restrict our attention only to odd $r$'s as the arguments are similar for even $r$'s.  By partial summation, we have
\begin{equation} \label{Marge}
\begin{split}
&\sum_{f \leq 2 \sqrt{x}} \frac{1}{f} 
\sum_{p \in S_{f,r}(r^2/(2\alpha)^2)} 
\frac{\sqrt{4p-r^2}}{p} L(1,\chi_d) \log p\\
 & = \frac{4\alpha\sqrt{1-\alpha^2}}{r} \sum_{f \leq 2 \sqrt{x}} \frac{1}{f} \sum_{p \in S_{f,r}(r^2/(2\alpha)^2)} L(1, \chi_d) \log p \\
& \hspace*{.25in} - \int\limits_{r^2/(2\beta)^2}^{r^2/(2\alpha)^2} \sum_{f \leq 2 \sqrt{x}} \frac{1}{f} \sum_{p \in S_{f,r}(t)} L(1, \chi_d) \log p \frac{\dif}{\dif t} \left(\frac{2}{\sqrt{t}}\sqrt{1-\frac{r^2}{4t}}\right) \dif t,
\end{split}
\end{equation}
where $d=(r^2-4p)/f^2$. \newline

For a fixed $U>0$ to be chosen later, we have
\begin{equation} \label{Lest}
L(1, \chi_d) = \sum_{n \leq U} \left( \frac{d}{n} \right) \frac{1}{n} + O \left( \frac{\sqrt{|d|} \log |d|}{U} \right) = \sum_{n \leq U} \left( \frac{d}{n} \right) \frac{1}{n} + O \left( \frac{\sqrt{|p|} \log |p|}{fU} \right),
\end{equation}
using the P\'olya-Vinogradov inequality, Lemma~\ref{PolyaVinolemma}. We further note that 
\begin{equation} \label{dist}
\frac{r^2}{(2\alpha)^2}-\frac{r^2}{(2\beta)^2} \ll \frac{r^2\gamma}{\beta\alpha^2}.
\end{equation}\newline

Using \eqref{Lest}, \eqref{dist} and the fact that $r \leq 2 \sqrt{x} \alpha$, we get that the double sum in the integrand on the right-hand side of \eqref{Marge} is
\begin{equation} \label{bart}
 \sum_{\substack{f \leq 2\sqrt{x} \\ (2r,f)=1}} \frac{1}{f} \sum_{n \leq U} \frac{1}{n} \sum_{p \in S_{f,r}(t)} \left( \frac{d}{n} \right) \log p + O \left( \frac{x^{3/2}\gamma \log^2 x}{\beta U} \right),
\end{equation}
where we have taken into account that $\gamma/\beta \geq 1/x$.
The coprimality relation $(2r,f)=1$ in the above comes in the following way.  As $r$ is odd, if $f^2|r^2-4p$, then $f$ is also odd and $d=(r^2-4p)/f^2 \equiv 1 \pmod{4}$.  Moreover since $(r,f)|p$, and $p \geq r$ if $p$ is neither 2 nor 3, we have $(r,f)=1$. The possible contribution of $p=2,3$ can be absorbed into 
the error term if $U\le x^{3/2}\gamma/\beta$ which will be the case by our later choice of $U$. \newline

For a fixed parameter $V$ with $1 \leq V \leq 2\sqrt{x}$ to be chosen later, the first term in \eqref{bart} is
\[  \sum_{\substack{f \leq V \\ (2r,f)=1}} \frac{1}{f} \sum_{n \leq U} \frac{1}{n} \sum_{p \in S_{f,r}(t)} \left( \frac{d}{n} \right) \log p +  \sum_{\substack{V< f \leq 2\sqrt{x} \\ (2r,f)=1}} \frac{1}{f} \sum_{n \leq U} \frac{1}{n} \sum_{p \in S_{f,r}(t)} \left( \frac{d}{n} \right) \log p. \]
Using \eqref{dist} and $r\le 2\sqrt{x}\alpha$, the second term of the above is 
estimated as
\begin{eqnarray*}
\left| \sum_{\substack{V< f \leq 2\sqrt{x} \\ (2r,f)=1}} \frac{1}{f} \sum_{n \leq U} \frac{1}{n} \sum_{p \in S_{f,r}(t)} \left( \frac{d}{n} \right) \log p \right| & \leq & \log x \log U \sum_{V<f \leq 2\sqrt{x}} \frac{1}{f} \sum_{\substack{r^2/(2\beta)^2 \leq n \leq t \\ n \equiv \bar{4}r^2 \bmod{f^2}}} 1 \\
& \ll & x \frac{\gamma}{\beta} \log x \log U \sum_{V<f \leq 2\sqrt{x}} \frac{1}{f^3} \\
& \ll & \frac{x \gamma \log x \log U}{\beta V^2},
\end{eqnarray*}
where $\bar{4}$ is an integer such that $\bar{4} \times 4 \equiv 1 \pmod{f^2}$. \newline

Therefore, \eqref{bart} may be re-written as
\begin{equation} \label{lisa}
\sum_{\substack{f \leq V \\ (2r,f)=1}} \frac{1}{f} \sum_{n \leq U} \frac{1}{n} \sum_{p \in S_{f,r}(t)} \left( \frac{d}{n} \right) \log p + O \left( \frac{\gamma}{\beta} x \log^2 x \left( \frac{x^{1/2}}{U} + \frac{1}{V^2} \right) \right),
\end{equation}
where we assume henceforth that $\log U \ll \log x$.  The sum over $f$ in \eqref{lisa} is evaluated by splitting the sum according to the residue class of $d$ modulo $4n$.  Since $d=(r^2-4p)/f^2$ is odd, and $\left( \frac{d}{n} \right)=0$ whenever $(d,n)>1$, we get
\begin{equation} \label{maggie}
 \sum_{\substack{f \leq V \\ (2r,f)=1}} \frac{1}{f} \sum_{n \leq U} \frac{1}{n} \sum_{p \in S_{f,r}(t)} \left( \frac{d}{n} \right) \log p = \sum_{\substack{n \leq U, f \leq V \\ (2r, f) =1}} \frac{1}{fn} \sum_{\substack{a \bmod{4n} \\ (a,4n)=1}} \left( \frac{a}{n} \right) \sum_{\substack{p \in S_{f,r}(t) \\ d \equiv a \bmod{4n}}} \log p. 
\end{equation}
The two conditions in the inner-most of the above sums $p \in S_{f,r}(t)$ and $d=(r^2-4p)/f^2 \equiv a \pmod{4n}$ are equivalent to $(r/2\beta)^2 \leq p \leq t$ and $p \equiv (r^2-af^2)/4 \pmod{nf^2}$.  Moreover, since $(2r,f)=1$, $nf^2$ and $(r^2-af^2)/4$ are co-prime if and only if $(r^2-af^2,4n)=4$. \newline

We re-write \eqref{maggie} as
\begin{equation} \label{abe}
\begin{split}
& \sum_{\substack{n \leq U, f \leq V \\ (2r, f) =1}} \frac{1}{fn} \sum_{\substack{a \bmod{4n} \\ (a,4n)=1}} \left( \frac{a}{n} \right) \Theta\left(r^2/(2\beta)^2 , t-r^2/(2\beta)^2; nf^2, \frac{r^2-af^2}{4} \right) \\
&= \left( t-\frac{r^2}{2\beta^2} \right) \sum_{\substack{n \leq U, f \leq V \\ (2r, f) =1}} \frac{c_f^r(n)}{fn \varphi(nf^2)} + \sum_{\substack{n \leq U, f \leq V \\ (2r, f) =1}} \frac{1}{fn} \\
& \hspace*{.5in} \times \sum_{\substack{a \bmod{4n} \\ (a,4n)=1 \\ (r^2-af^2,4n)=4}} \left( \frac{a}{n} \right)  E\left(r^2/(2\beta)^2 , t-r^2/(2\beta)^2; nf^2, \frac{r^2-af^2}{4} \right),
\end{split}
\end{equation}
where $c_f^r(n)$, $\Theta(x,y;q,a)$ and $E(x,y;q,a)$ are as defined in
 \eqref{constantfrn} and \eqref{defpsi1}. \newline

We first deal with the case 
\begin{equation} \label{milhouse}
 t- \left( \frac{r}{2 \beta} \right)^2 \geq \left( \frac{r}{2 \beta} \right)^{7/6} x^{\varepsilon}.
\end{equation}
In this case, applying Cauchy's inequality, we have that the second term in \eqref{abe} is
\begin{eqnarray*}
& \ll & \sum_{\substack{f \leq V \\ (2r,f)=1}} \frac{1}{f} \left( \sum_{n \leq U} \frac{\varphi(4n)}{n^2} \right)^{\frac{1}{2}} \left( \sum_{n \leq U} \sum_{\substack{a \bmod{4n} \\ (a,4n)=1 \\ (r^2-af^2,4n)=4}} E^2 \left( \frac{r^2}{(2\beta)^2} , t-\frac{r^2}{(2\beta)^2}; nf^2, \frac{r^2-af^2}{4} \right) \right)^{\frac{1}{2}} \\
& \ll & ( \log U)^{\frac{1}{2}} \sum_{\substack{f \leq V \\ (2r,f)=1}} \frac{1}{f} \left( \sum_{n \leq U} \sum_{\substack{b \bmod{nf^2} \\ (b,nf^2)=1}} E^2 \left( \frac{r^2}{(2\beta)^2} , t-\frac{r^2}{(2\beta)^2}; nf^2, b \right) \right)^{\frac{1}{2}},
\end{eqnarray*}
as $a_1 \not\equiv a_2 \pmod{4n}$ ensures that 
\[ b_1 = (r^2-a_1f^2)/4 \not\equiv b_2=(r^2-a_2f^2)/4 \pmod{nf^2}. \]

The last sum above is majorized by
\[ \log V (\log U)^{\frac{1}{2}} \left( \sum_{n \leq UV^2} \sum_{\substack{a \bmod{n} \\ (a,n)=1}} E^2 \left(\frac{r^2}{(2\beta)^2} , t-\frac{r^2}{(2\beta)^2}; n, a \right) \right)^{\frac{1}{2}}. \]
Fix any $c>0$.  Because of the condition in \eqref{milhouse}, we may apply Lemma \ref{BDVlemma} to the double sum above and it is bounded by
\begin{equation} \label{nelson}
 \ll \frac{x \gamma}{\beta \log^c x}
\end{equation}
if 
\begin{equation} \label{martin}
 \log U \ll \log x, \; \log V \ll \log x, \; r \leq 2 \sqrt{x} \alpha \; \mbox{and} \; UV^2 \leq \frac{x \gamma}{\beta (\log x)^{5+2c}}.
\end{equation}
Here we have again used \eqref{dist}. \newline

If \eqref{milhouse} is not satisfied, then the expression of interest satisfies the same majorant given in \eqref{nelson} by a trivial bound for 
\[ E^2 \left( r^2/(2\beta)^2 , t-r^2/(2\beta)^2; nf^2, \frac{r^2-af^2}{4} \right), \]
where we take into account the estimate (\ref{dist}) and the conditions  
$\gamma/\beta \le x^{-\varepsilon}$ and   
\[ \frac{x^{2\varepsilon}\alpha \beta^{6/5}}{\gamma^{6/5}}<r\le 
2\sqrt{x}\alpha. \]
We now set, assuming $c>1$,
\[ U = x^{1/2} \log^{c+2} x, \; V = \log^c x. \]
Note that these choices for $U$ and $V$ are consistent with the conditions in \eqref{martin}, since we assume that
\[ \frac{\gamma}{\beta} \geq x^{-1/2+\varepsilon}. \]
Combining everything and using Lemma \ref{maintermlemma}, the second double sum on the right-hand side of \eqref{Marge} can be expressed in the following way.
\[
\sum_{f \leq 2 \sqrt{x}} \frac{1}{f} \sum_{p \in S_{f,r}(t)} L(1, \chi_d) \log p = \left( t- \frac{r^2}{2\beta^2} \right) K_r + O \left( \frac{x \gamma}{\beta \log^c x} \right). 
\]
Therefore, the left-hand side of \eqref{Marge} is
\begin{equation} \label{theend}
=  \frac{\pi K_rrF(\alpha,\beta)}{\alpha^2}\left(1+O\left(
x^{-\varepsilon}\right)\right)+O\left(\frac{xF(\alpha,\beta)}{r\log^c x}
\right),
\end{equation}
where for the error term estimate we have used
\[ \gamma \sqrt{1-\alpha^2} \ll F(\alpha, \beta), \]
and for the main term we have applied integration by parts, the change of variables
$v^2=r^2/4t$ and 
\[ \frac{r}{v^2} = \frac{r}{\alpha^2} \left( 1 + O(x^{-\varepsilon}) \right), \; \mbox{for} \; \alpha \leq v \leq \beta, \]
which follows from \eqref{dist} and $\gamma/\beta \leq x^{-\varepsilon}$.  Summing up \eqref{theend} over $r$ in the interval 
$x^{2\varepsilon}\alpha\beta^{6/5}/ \gamma^{6/5} < r \le 2\sqrt{x}\alpha$, 
using Lemma \ref{LTconstantlemma}
and partial summation, and taking into account 
the conditions $\gamma/\beta\ge
x^{2\varepsilon-5/12}$ and $\alpha\ge x^{-1/2+\varepsilon}$, 
\eqref{Burns} and hence $S_2$ are equal to
\begin{equation} \label{seconddoublesum}
xF(\alpha,\beta)\left(1+\frac{1}{\log^{c-1} x}\right).
\end{equation}  

It remains to show that $S_1$ and $S_3$ are small comparing to $S_2$.  By Lemma \ref{Hlemma} and the well-known estimate $L(1,\chi_d)\ll \log 2|d|$, we have
\begin{equation}\label{Heq}
H(r^2-4p)\ll \sqrt{|r^2-4p|}\log^2 x
\end{equation}
from which we deduce, by a short calculation, that $S_1$ and $S_3$ are majorized by
\begin{equation} \label{firstds}
\ll x^{1-\varepsilon}F(\alpha,\beta).
\end{equation}
Here we have used the condition $\gamma/\beta\ge x^{2\varepsilon-5/12}$ for the estimation of $S_1$ and the condition $\gamma/\beta\le x^{-\varepsilon}$ for the estimation of $S_3$.  Now combining \eqref{Homer}, \eqref{seconddoublesum} and \eqref{firstds}, we obtain the desired result.
\end{proof}
 
\section{Preparations for the estimation of the error term}
We first characterize the elliptic curves lying in a fixed $\fie_p$-isomorphism class, where $p$ is a prime that is neither 2 nor 3.  In the following, for $z\in \intz$ let $\overline{z}$ be the reduction of $z$ modulo $p$. Furthermore, let $z^{-1}$ be a multiplicative inverse modulo $p$, that is, $z z^{-1}\equiv 1 \pmod{p}$. \\

\begin{lemma} \label{charlemma}
Let $a,b,c,d\in \intz$, $p\nmid abcd$ and $E_1$,
$E_2$ be  
elliptic curves over $\fie_p$ given by
$$
E_1\ :\ Y^2=X^3+\overline{a}X+\overline{b} \; \mbox{and} \; 
E_2\ :\ Y^2=X^3+\overline{c}X+\overline{d}.
$$

{\bf (i)} If $p\equiv 1 \pmod{4}$, then $E_1$ and $E_2$ are $\fie_p$-isomorphic if and only if $ca^{-1}$ is a biquadratic residue modulo $p$ and $c^3a^{-3}\equiv d^2b^{-2} \pmod{p}$.\medskip

{\bf (ii):} If $p\equiv 3 \pmod{4}$, then $E_1$ and $E_2$ are $\fie_p$-isomorphic if and only if $ca^{-1}$ and $db^{-1}$ are quadratic residues modulo $p$ and $c^3a^{-3}\equiv d^2b^{-2} \pmod{p}$.
\end{lemma}
\begin{proof} By Lemma \ref{LE2}, the curves $E_1$ and $E_2$ are $\fie_p$-isomorphic if and only if there exists an integer $m$ such that $p\nmid m$ and 
\begin{equation} \label{bed}
c\equiv m^4a \pmod{p} \ \ \ \ \mbox{and } \ \ \ d\equiv m^6b \pmod{p}.
\end{equation}

{\bf (i)} Suppose that $p\equiv 1 \pmod{4}$.  If \eqref{bed} is satisfied, then it follows that $ca^{-1}$ is a biquadratic residue modulo $p$ and $c^3a^{-3}\equiv m^{12}\equiv d^2b^{-2} \pmod{p}$. \newline

Assume, conversely, that  $ca^{-1}$ is a biquadratic residue modulo $p$ and 
\begin{equation} \label{congru}
c^3a^{-3}\equiv d^2b^{-2} \pmod{p}. 
\end{equation}

Since $p\equiv 1 \pmod{4}$, there exist two solutions $m_1,m_2$ of the congruence $c\equiv m^4a \pmod{p}$ such that $m_2^2\equiv -m_1^2 \pmod{p}$, and \eqref{congru} implies that $d^2b^{-2}\equiv m_j^{12} \pmod{p}$ for $j=1,2$. From this it follows that $db^{-1}\equiv m_1^6 \pmod{p}$ or $db^{-1}\equiv -m_1^6\equiv 
m_2^6 \pmod{p}$.  Hence, the system \eqref{bed} is solvable for $m$. This completes the proof of {\bf (i)}.
\medskip

{\bf (ii)} Suppose that $p\equiv 3 \pmod{4}$.  If (\ref{bed}) is satisfied, then it follows that $ca^{-1}$ and $db^{-1}$ are quadratic residues modulo $p$ and $c^3a^{-3}\equiv m^{12}\equiv d^2b^{-2} \pmod{p}$. \newline

Assume, conversely, that $ca^{-1}$ and $db^{-1}$ are quadratic residues modulo $p$ and (\ref{congru}) is satisfied. Then,since $p\equiv 3 \pmod{4}$, $ca^{-1}$ is also a biqadratic residue.  Hence, there exists a solution $m$ of the congruence $c\equiv m^4a \pmod{p}$. Further, \eqref{congru} implies that $d^2b^{-2}\equiv m^{12} \pmod{p}$. From this it follows that $db^{-1}\equiv m^6 \pmod{p}$ or $db^{-1}\equiv -m^6 \pmod{p}$. But $-m^6$ is a quadratic non-residue modulo $p$ since $p\equiv 3 \pmod{4}$. Thus $db^{-1}\not\equiv -m^6 \pmod{p}$ since $db^{-1}$ is supposed to be a quadratic residue modulo $p$. Hence, we have $db^{-1}\equiv m^6 \pmod{p}$, and so the system \eqref{bed} is solvable for $m$. This completes the proof of {\bf (ii)}. \end{proof}

We shall detect elliptic curves lying in a fixed $\fie_p$-isomorphism class by using Dirichlet characters. For the estimation of certain error terms we then need the following results on character sums.

\begin{lemma} \label{Dirlemma} 
Let $q,N\in \mathbbm{N}$ and $(a_n)$ be any sequence of complex numbers. Then 
$$
\sum\limits_{\chi\ \! \bmod{q}} \left| \sum\limits_{n\le N} a_n\chi(n) \right|^2
=\varphi(q)\sum\limits_{\substack{a=1\\ (a,q)=1}}^{q} \left| \sum\limits_{\substack{n\le N\\ n\equiv a\ \!\bmod{q}}} a_n\right|^2,
$$ 
where the outer 
sum on the left-hand side runs over all Dirichlet characters modulo $q$.
\end{lemma}
\begin{proof} This is a consequence of the orthogonality relations for 
Dirichlet characters. \end{proof}

\begin{lemma}\label{fourthmomentlemma} 
Let $q,N\in \mathbbm{N}$, $q\ge 2$. Then 
$$
\sum\limits_{\chi\not=\chi_0} \left| \sum\limits_{n\le N} \chi(n) \right|^4
\ll N^2q\log^6 q,
$$
where the outer sum on the left-hand side runs over all non-principal 
Dirichlet characters
modulo $q$.
\end{lemma}
\begin{proof} This is Lemma 3 in \cite{FrIw}. \end{proof}

Furthermore, we shall need the following estimates for sums over 
$H_{r,p}$.\\

\begin{lemma} \label{HPlemma} Suppose that $0<\alpha\le \beta\le 1$,
$x\ge 1$ and $F(\alpha,\beta)\ge x^{-1/2}$, where $F(\alpha,\beta)$ is
defined as in \eqref{Fdef}. Let 
$$
H_p:=\sum\limits_{2\sqrt{p}\alpha\le r \le 2\sqrt{p}\beta}
H(r^2-4p).
$$
Then
$$
\sum\limits_{p\le x} H_{p}^{1/2}\ll x^{3/2+\varepsilon}F(\alpha,\beta)^{1/2}, 
\ \ \ \ \ \
\sum\limits_{p\le x} \frac{H_{p}}{\sqrt{p}}\ll 
x^{3/2+\varepsilon}F(\alpha,\beta),
$$ 
and 
$$
\sum\limits_{p\le x} \frac{H_{p}}{p}\ll x^{1+\varepsilon}F(\alpha,\beta),
 \ \ \ \ \ \ \sum\limits_{p\le x} \frac{H_{p}\log p}{p^2}\ll 1.
$$
\end{lemma}

\begin{proof} By Lemma \ref{Hlemma}, we have 
\begin{equation} \label{convert}
H_p
=\sum\limits_{2\sqrt{p}\alpha\le r \le 2\sqrt{p}\beta}
\frac{1}{\pi}\sum\limits_{\substack{f,d\\ r^2-4p=df^2\\ d\equiv 0,1\
\bmod{4}}} \sqrt{|d|}L(1,\chi_d).
\end{equation}
It is well-known that $L(1,\chi_d)\ll \log 2|d|$ and that the number of divisors $\tau(n)$ of a natural number $n$ satisfies $\tau(n)\ll n^{\varepsilon}$ . Therefore, from \eqref{convert} it follows that
\begin{equation} \label{Hest}
H_p \ll \sum\limits_{2\sqrt{p}\alpha\le r \le 2\sqrt{p}\beta} (4p-r^2)^{1/2+\varepsilon} \ll p^{1/2+\varepsilon} \sum\limits_{2\sqrt{p}\alpha\le r \le 2\sqrt{p}\beta} \sqrt{1-\left(\frac{r}{2\sqrt{p}}\right)^2} \ll p^{1/2+\varepsilon} \left(1+\sqrt{p}F(\alpha,\beta)\right).
\end{equation}
From \eqref{Hest}, we obtain
\begin{equation}\label{HE}
\sum\limits_{p\le x} H_{p} \ll x^{2+\varepsilon}F(\alpha,\beta)
\end{equation}
if $F(\alpha,\beta)\ge x^{-1/2}$. 
Using the Cauchy-Schwarz inequality, we obtain
$$
\sum\limits_{p\le x} H_{p}^{1/2}\ll 
x^{1/2}\left(\sum\limits_{p\le x} H_{p}\right)^{1/2}\ll  x^{3/2+\varepsilon}
F(\alpha,\beta)^{1/2}
$$
from \eqref{HE}. The remaining three estimates in Lemma 13 can be derived from \eqref{HE} by partial summation.
\end{proof}

Finally, we shall need the following bound.

\begin{lemma} \label{numb} 
The number of $\fie_p$-isomorphism classes of elliptic curves 
containing curves  
$$
E\ :\ Y^2=X^3+aX+b
$$
over $\fie_p$  with $a=0$ or $b=0$ is bounded by 10.
\end{lemma}

\begin{proof} By  Lemma 2, the number of $\fie_p$-isomorphism classes containing curves $E(0,b)$ with $b\in \fie_p^* $ is less than 6, and the number of $\fie_p$-isomorphism classes containing curves $E(a,0)$ with $a\in \fie_p^*$ does not exceed 4. \end{proof} 
   
\section{Estimation of the error term}
Now we evaluate the error term more explicitly than in \eqref{DaPa}. We shall
establish the following.
\begin{theorem} \label{errortermtheorem}
Let $A,B,x\ge 1$ and $0<\alpha\le \beta\le 1$ Suppose that 
$F(\alpha,\beta)\ge x^{-1/2}$, where $F(\alpha,\beta)$ is
defined as in \eqref{Fdef}. Then
\begin{equation*}
\begin{split}
\frac{1}{4AB}\sum\limits_{|a|\le A} &\sum\limits_{|b|\le B}
\Theta_{E(a,b)}(\alpha,\beta;x) -
\sum\limits_{p\le x}\log p \sum\limits_{2\sqrt{p}\alpha\le 
r\le 2\sqrt{p}\beta} \frac{H(r^2-4p)}{2p}\\
&\ll x^{3/2+\varepsilon}\left(F(\alpha,\beta)
\left(A^{-1}+B^{-1}\right)+ F(\alpha,\beta)^{1/2}(AB)^{-1/2}\right).
\end{split}
\end{equation*}
\end{theorem}

\begin{proof}
For $p>3$ let $I_{p}$ be the number of $\fie_p$-isomorphism classes of elliptic curves 
$$
E\ :\ Y^2=X^3+cX+d
$$
over $\fie_p$ with $2\sqrt{p}\alpha\le \lambda_E(p)\le 2\sqrt{p}\beta$ points
such that $c,d\not=0$. Let $(u_{p,j},v_{p,j})$, $j=1,...,I_{p}$ be pairs of
integers such that the curves $E(\overline{u_{p,j}},\overline{v_{p,j}})$ 
form a 
system of representatives of these isomorphism classes. 
We now write
\begin{eqnarray*}
& & \sharp\{|a|\le A,\ |b|\le B\ :\ 2\sqrt{p}\alpha\le 
\lambda_{E(a,b)}(p)\le 2\sqrt{p}\beta\}\\  &=&
\sharp\{|a|\le A,\ |b|\le B\ :\ p \nmid ab,\  2\sqrt{p}\alpha\le 
\lambda_{E(a,b)}(p)\le 2\sqrt{p}\beta\} + O\left(\frac{AB}{p}+A+B\right)\nonumber
\end{eqnarray*}
and 
\begin{equation}\label{otto}
\sharp\{|a|\le A, \ |b|\le B\ :\ p\nmid ab,\  2\sqrt{p}\alpha\le \lambda_{E(a,b)}(p)\le 2\sqrt{p}\beta\}
= \sum\limits_{j=1}^{I_p} \sharp\{|a|\le A,\ |b|\le B\ : E(\overline{a},\overline{b})\cong 
E(\overline{u_{p,j}},\overline{v_{p,j}})\},
\end{equation}
where  the symbol $\cong$ stands for ``$\fie_p$-isomorphic''.  We rewrite the term on the right-hand  side of (\ref{otto}) as a character sum.  If $p\equiv 1 \pmod{4}$, then, by Lemma~\ref{charlemma} (i) and the character relations, 
this term equals
\begin{equation} \label{character}
\frac{1}{4\varphi(p)}\sum\limits_{j=1}^{I_{p}} 
\sum\limits_{|a|\le A} \sum\limits_{|b|\le B} \sum\limits_{k=1}^4 
\left(\frac{au_{p,j}^{-1}}{p}\right)_4^k \
\sum\limits_{\chi\ \! \bmod{p}} 
\chi(a^3u_{p,j}^{-3}b^{-2}v_{p,j}^2),
\end{equation}
where $(\cdot/p)_4$ is the biquadratic residue symbol.  If $p\equiv 3 \pmod{4}$, then, by Lemma~\ref{charlemma} (ii) and the character relations, the term on the right-hand  side of \eqref{otto} equals
\[ \frac{1}{4\varphi(p)}\sum\limits_{j=1}^{I_{p}} 
\sum\limits_{|a|\le A} \sum\limits_{|b|\le B} \left(\chi_0(a)+
\left(\frac{au_{p,j}^{-1}}{p}\right)\right) 
\left(\chi_0(b)+
\left(\frac{bv_{p,j}^{-1}}{p}\right)\right) \sum\limits_{\chi\ \! \bmod{p}} 
\chi(a^3u_{p,j}^{-3}b^{-2}v_{p,j}^2), \]
where $(\cdot/p)$ is the Legendre symbol and $\chi_0$ is the principal character.\newline

In the following, we consider only the case $p\equiv 1 \pmod{4}$. The case $p\equiv 3 \pmod{4}$ can be treated in a similar way.  The expression in \eqref{character} equals
$$
\frac{1}{4\varphi(p)} \sum\limits_{k=1}^4 \sum\limits_{\chi\ \! \bmod{p}} \sum\limits_{j=1}^{I_{p}} \left(\frac{u_{p,j}}{p}\right)_4^{-k} \overline{\chi}^3(u_{p,j}) \chi^2(v_{p,j}) \sum\limits_{|a|\le A}\left(\frac{a}{p}\right)_4^{k}\chi^3(a) \sum\limits_{|b|\le B}  \overline{\chi}^2(b).
$$
We split this expression into 3 parts $M,E_1,E_2$, where\\ \\
({\bf i}) $M=$ contribution of $k,\chi$ with 
$(\cdot/p)_4^{k} \chi^3=\chi_0$, 
$\chi^2=\chi_0$;\\ \\
({\bf ii}) $E_1=$ contribution of $k,\chi$ with 
$(\cdot/p)_4^{k} \chi^3\not=\chi_0$, 
$\chi^2=\chi_0$ or 

\ $(\cdot/p)_4^{k} \chi^3=\chi_0$, 
$\chi^2\not=\chi_0$;\\ \\
({\bf iii}) $E_2=$ contribution of $k,\chi$ with 
$(\cdot/p)_4^{k} \chi^3\not=\chi_0$, $\chi^2\not=\chi_0$.\\ \\
As one may expect, $M$ shall turn out to be the main term and $E_1$, $E_2$
to be the error terms.\\

{\it Estimation of $M$.}\ \ \ \
The only cases in which 
$(\cdot/p)_4^{k} \chi^3=\chi_0 \mbox{ and } \chi^2=\chi_0$ are 
$k=0$, $\chi=\chi_0$ and $k=2$, $\chi=(\cdot/p)$. Now, by a short calculation,
we obtain
\begin{equation} \label{S1}
M=\frac{4ABI_{p}}{2p}\left(1+O\left(\frac{1}{p}\right)\right).
\end{equation}
By Lemma \ref{numb}, we have 
$$
\sum\limits_{2\sqrt{p}\alpha\le r \le 2\sqrt{p}\beta}
H(r^2-4p)-I_{p}\le 10.
$$
Combining this with \eqref{S1}, we obtain
$$
M=\sum\limits_{2\sqrt{p}\alpha\le r\le 2\sqrt{p}\beta} 
\frac{4ABH(r^2-4p)}{2p}+ O\left(\frac{AB}{p}+\frac{ABI_{p}}{p^2}\right).
$$\medskip

{\it Estimation of $E_1$.}\ \ \ \ The number of solutions $(k,\chi)$ with $k=1,...,4$ of $(\cdot/p)_4^{k} \chi^3=\chi_0$ is bounded by 12, and $\chi^2=\chi_0$ has precisely 2 solutions $\chi$.  Thus $E_1$ is the sum of at most $12+4\times 2=20$ terms of the form
$$ 
\frac{1}{4\varphi(p)} \sum\limits_{j=1}^{I_{p}} 
\overline{\chi_1}(u_{p,j})\overline{\chi_2}(v_{p,j})
\sum\limits_{|a|\le A} \chi_1(a) \sum\limits_{|b|\le B} \chi_2(b),
$$
where exactly one of the characters $\chi_1$, $\chi_2$ is the principal character $\chi_0$. Therefore, Lemma~\ref{PolyaVinolemma} implies that
$$
E_1\ll \frac{I_{p}(A+B)}{\sqrt{p}}\log p.
$$ 
 
{\it Estimation of $E_2$.}\ \ \ \  Given $k\in \intz$ and a character $\chi_1 \pmod{p}$, the number of solutions $\chi$ of $\left(\frac{\cdot}{p}\right)_4^{k} {\chi}^{-3}=\chi_1$ is less than 3, and the number of solutions $\chi$ of $\chi^2=\chi_1$ is $\le 2$.  Thus, using the Cauchy-Schwarz inequality, we deduce that
\begin{equation} \label{cha}
E_2 \ll \frac{1}{p} \sum\limits_{k=1}^4 \left( \sum\limits_{\chi} \left|\sum\limits_{j=1}^{I_{p}} \left(\frac{u_{p,j}}{p}\right)_4^{k} \chi(u_{p,j}^{-3}v_{p,j}^2) \right|^2 \right)^{1/2} \left(\sum\limits_{\chi\not=\chi_0} \left| \sum\limits_{|a|\le A} \chi(a) \right|^4\right)^{1/4} \left(\sum\limits_{\chi\not=\chi_0} \left| \sum\limits_{|b|\le B}  \chi(b)\right|^4\right)^{1/4}.
\end{equation}  
By Lemma~\ref{charlemma} (i), the number of $j$'s such that $u_{p,j}^{-3}v_{p,j}^2$ lie in a fixed residue class modulo $p$ is bounded by 4.  Using this, Lemma~\ref{Dirlemma} and Lemma~\ref{fourthmomentlemma}, the expression on the right-hand side of \eqref{cha} is dominated by 
$$
\ll (I_{p}AB)^{1/2}\log^3 p.
$$

{\it The final estimate.} Combining all contributions, we obtain 
\begin{equation} \label{end}
\begin{split}
& \sharp\{|a|\le A, \ |b|\le B\ :\ 2\sqrt{p}\alpha\le 
\lambda_{E(a,b)}(p)\le 2\sqrt{p}\beta\} \\ 
&= \sum\limits_{2\sqrt{p}\alpha\le r\le 2\sqrt{p}\beta} 
\frac{4ABH(r^2-4p)}{2p} +
O\left(\frac{AB}{p}+ \frac{ABI_{p}}{p^2}+A+B+
\left(I_{p}AB\right)^{\frac{1}{2}}\log^3 p+ \frac{I_{p}(A+B)}{\sqrt{p}}\log p\right).
\end{split}
\end{equation}
The result of Theorem~\ref{TH1} now follows from \eqref{goal}, \eqref{end}, Lemma \ref{HPlemma}, $I_p\le H_p$ and the prime number theorem which says that $\pi(x)\sim x/\log x$ as $x\rightarrow \infty$. \end{proof}

Combining Theorem \ref{maintermtheorem} and Theorem \ref{errortermtheorem},
we obtain Theorem \ref{TH1}. Note that the condition 
$\alpha\ge x^{-1/2+\varepsilon}$ in Theorem
\ref{maintermtheorem} can be dropped if we 
assume the conditions 
$x^{\varepsilon-5/12}\le \gamma/\beta \le x^{-\varepsilon}$ and $F(\alpha,
\beta)\ge x^{-1/2+\varepsilon}$ in Theorem \ref{TH1}.

\section{Proof of Theorem 2}
We set
$$
\mu:=\frac{1}{4AB}\sum\limits_{|a|\le A} \sum\limits_{|b|\le B}
\Theta_{E(a,b)}(\alpha,\beta;x).
$$
Fix any $c>0$ and assume that $A,B>x^{1/2+\varepsilon}$ and 
$AB>x^{1+\varepsilon}/F(\alpha,\beta)$. Then, by Corollary 1,
we have 
\begin{equation}\label{mu}
\mu = xF(\alpha,\beta)+O\left(x\frac{F(\alpha,\beta)}{\log^c x}\right).
\end{equation}
Thus, by the triangle inequality, the left-hand side of (\ref{T3}) is
\begin{equation} \label{start}
\begin{split}
& \ll \frac{1}{4AB}\left|\sum\limits_{|a|\le A} \sum\limits_{|b|\le B}
\Theta_{E(a,b)}(\alpha,\beta;x)
-\mu\right|^2+O\left(\frac{x^2F(\alpha,\beta)^2}{\log^{2c} x}\right)\\ 
&= \frac{1}{4AB}\sum\limits_{|a|\le A} \sum\limits_{|b|\le B}
\Theta_{E(a,b)}(\alpha,\beta;x)^2
-\mu^2+O\left(\frac{x^2F(\alpha,\beta)^2}{\log^{2c} x}\right),
\end{split}
\end{equation}
where the second line arises from the general identity
$$
\frac{1}{N}\sum\limits_{n=1}^N \left(a_n-\mu\right)^2=
\frac{1}{N}\sum\limits_{n=1}^N a_n^2-\mu^2\ \ \ \ \ \mbox{ if }\ \ \ \ \
\mu=\frac{1}{N} \sum\limits_{n=1}^N a_n.
$$
We now write 
\begin{equation} \label{goal2}
\sum\limits_{|a|\le A} \sum\limits_{|b|\le B} \Theta_{E(a,b)}(\alpha,\beta;x)^2= 
4AB\mu + \sum\limits_{\substack{p,q\le  x\\ p\not= q}} \log p \log q \mathop{\sum \sum}_{\substack{|a|\le A,\ |b|\le B \\ 2\sqrt{p}\alpha\le \lambda_{E(a,b)}(p)\le 2\sqrt{p}\beta \\  2\sqrt{q}\alpha\le \lambda_{E(a,b)}(q)\le 2\sqrt{q}\beta }} 1,
\end{equation}
where $p,q$ denote primes.  Obviously, 
\begin{equation} \label{app}
\begin{split}
\sum\limits_{\substack{p,q\le  x\\ p\not= q}} & \log p\log q \mathop{\sum \sum}_{\substack{|a|\le A,\ |b|\le B \\ 2\sqrt{p}\alpha\le \lambda_{E(a,b)}(p)\le 2\sqrt{p}\beta \\  2\sqrt{q}\alpha\le \lambda_{E(a,b)}(q)\le 2\sqrt{q}\beta }} 1 \\
&= \sum\limits_{\substack{p,q\le  x\\ p\not= q}} \log p \log q 
\mathop{\sum \sum}_{\substack{|a|\le A,\ |b|\le B,\ p, q \nmid ab \\ 2\sqrt{p}\alpha\le \lambda_{E(a,b)}(p)\le 2\sqrt{p}\beta \\  2\sqrt{q}\alpha\le \lambda_{E(a,b)}(q)\le 2\sqrt{q}\beta }} 1 +O\left( \sum\limits_{p\le x} \ \sum\limits_{\substack{|a|\le A,\ \! |b|\le B\\
p|ab}} \Theta_{E(a,b)}(\alpha,\beta;x)\right).
\end{split}
\end{equation}
Using Corollary 1 and $\sharp\{p\ :\ p|ab\}=\omega(|ab|)\ll \log\log(10|ab|)$ 
if $ab\not=0$, we have
\begin{equation} \label{app1}
\sum\limits_{p\le x} \ \sum\limits_{\substack{|a|\le A,\ \! |b|\le B\\p|ab}} \Theta_{E(a,b)}(\alpha,\beta;x) + O\left(xABF(\alpha,\beta)\log\log(10AB) +x^2(A+B)\right).
\end{equation}

Now we fix $p,q$ with $p\not=q$. 
In the following, we confine ourselves to the case when $p\equiv q\equiv 1 \pmod{4}$. The remaining cases $pq\equiv -1 \pmod{4}$ and $p\equiv q\equiv 3 \pmod{4}$ can be treated in a similar way. 
Similarly as in section 6, we can express the term
\[ \mathop{\sum \sum}_{\substack{|a|\le A,\ |b|\le B,\ p, q \nmid ab \\ 2\sqrt{p}\alpha\le \lambda_{E(a,b)}(p)\le 2\sqrt{p}\beta \\  2\sqrt{q}\alpha\le \lambda_{E(a,b)}(q)\le 2\sqrt{q}\beta }} 1 \]
as a character sum
\[ \frac{1}{16\varphi(p)\varphi(q)}
\sum\limits_{i=1}^{I_{p}}\sum\limits_{j=1}^{I_{q}} 
\sum\limits_{|a|\le A} \sum\limits_{|b|\le B} \sum\limits_{k=1}^4
\left(\frac{au_{p,i}^{-1}}{p}\right)_4^k\
\sum\limits_{\chi\ \! \bmod{p}} 
\chi(a^3u_{p,i}^{-3}b^{-2}v_{p,i}^2) \sum\limits_{l=1}^4 
\left(\frac{au_{q,j}^{-1}}{q}\right)_4^l\
\sum\limits_{\chi'\ \! \bmod{q}} 
\chi'(a^3u_{q,j}^{-3}b^{-2}v_{q,j}^2).\]
This sum equals
\begin{equation} \label{character3}
\begin{split}
& \frac{1}{16\varphi(p)\varphi(q)}\sum\limits_{k=1}^4 \sum\limits_{l=1}^4
\ \sum\limits_{\chi\ \! \bmod{p}}\
\sum\limits_{\chi'\ \! \bmod{q}} 
\left(\sum\limits_{i=1}^{I_{p}}
\left(\frac{u_{p,i}}{p}\right)_4^{-k}
{\overline\chi}^3(u_{p,i})\chi^2(v_{p,i})\right) \\
& \times \left(\sum\limits_{j=1}^{I_{q}}
\left(\frac{u_{q,j}}{q}\right)_4^{-l}
\overline{\chi'}^3(u_{q,j}){\chi'}^2(v_{q,j})\right)
\left(\sum\limits_{|a|\le A} 
\left(\frac{a}{p}\right)_4^k\left(\frac{a}{q}\right)_4^l 
\left(\chi{\chi'}\right)^3(a)\right) \left(\sum\limits_{|b|\le B} 
\left(\overline{\chi\chi'}\right)^2(b)\right).
\end{split}
\end{equation}
Let $\chi_0$ be the principal character modulo $p$ and 
$\chi'_0$ be the principal character modulo $q$. Then $\chi_0\chi'_0$ is
the principal character modulo $pq$. As in section 6, we split the
expression in (\ref{character3}) into 3 parts $M,E_1,E_2$, where\\ \\
({\bf i}) $M=$ contribution of $k,l,\chi,\chi'$ with 

\[ (\cdot/p)_4^{k} (\cdot/q)_4^{l}(\chi\chi')^3=\chi_0\chi'_0, (\chi\chi')^2=\chi_0\chi'_0; \]
({\bf ii}) $E_1=$ contribution of $k,l,\chi,\chi'$ with 

\[ (\cdot/p)_4^{k} (\cdot/q)_4^{l}(\chi\chi')^3\not=\chi_0\chi'_0, (\chi\chi')^2=\chi_0\chi'_0 \; \mbox{or} (\cdot/p)_4^{k} (\cdot/q)_4^{l}(\chi\chi')^3=\chi_0\chi'_0, (\chi\chi')^2\not=\chi_0\chi'_0; \]
({\bf iii}) $E_2=$ contribution of $k,l,\chi,\chi'$ with 

\[ (\cdot/p)_4^{k}(\cdot/q)_4^{l}(\chi\chi')^3\not=\chi_0\chi'_0, (\chi\chi')^2\not=\chi_0\chi'_0. \]
\\ 

{\it Estimation of $M$.}\ \ \ \
The only cases in which 
$(\cdot/p)_4^{k} (\cdot/q)_4^{l}(\chi\chi')^3=\chi_0\chi'_0$, 
$(\chi\chi')^2=\chi_0\chi'_0$ are:\medskip\\
(a) $k=l=0$, $\chi=\chi_0$, $\chi'=\chi'_0$; \medskip\\
(b) $k=l=2$, $\chi=(\cdot/p)$, $\chi'=(\cdot/q)$;\medskip\\
(c) $k=0$, $l=2$, $\chi=\chi_0$, $\chi'=(\cdot/q)$;\medskip\\
(d) $k=2$, $l=0$, $\chi=(\cdot/p)$, $\chi'=\chi_0$.\medskip\\ 
Now, by a short calculation, we obtain
\begin{equation} \label{S2}
M=\frac{4ABI_{p}I_{q}}{4pq}\left(1+O\left(\frac{1}{p}+\frac{1}{q}\right)
\right).
\end{equation}
By Lemma \ref{numb}, we have 
$$
\sum\limits_{2\sqrt{p}\alpha\le r \le 2\sqrt{p}\beta}
H(r^2-4p)-I_{p}\le 10 \ \ \ \ \ \ \mbox{and} \ \ \ \ \ \ 
\sum\limits_{2\sqrt{q}\alpha\le r \le 2\sqrt{q}\beta}
H(r^2-4q)-I_{q}\le 10.
$$
Combining this with (\ref{S2}), we obtain
\[ M=4AB\sum\limits_{2\sqrt{p}\alpha\le r \le 2\sqrt{p}\beta}
\frac{H(r^2-4p)}{2p}\sum\limits_{2\sqrt{q}\alpha\le r \le 2\sqrt{q}\beta}
\frac{H(r^2-4q)}{2q} +O\left(\frac{AB(I_{p}+I_{q})}{pq}
+ABI_{p}I_{q}\left(\frac{1}{p^2q}+\frac{1}{pq^2}\right)\right). \]

{\it Estimation of $E_1$.}\ \ \ \ The number of solutions 
$(k,l,\chi,\chi')$ with
$k,l=1,...,4$ of 
$(\cdot/p)_4^{k} (\cdot/q)_4^{l}(\chi\chi')^3\not=\chi_0\chi'_0$
is bounded by $12^2$, and $(\chi\chi')^2=\chi_0\chi_0'$
has precisely 4 solutions $(\chi,\chi')$.
Thus $E_1$ is the sum of at most $144+16\cdot 4=228$
terms of the form
$$ 
\frac{1}{16\varphi(p)\varphi(q)} \sum\limits_{|a|\le A} \chi_1(a) 
\sum\limits_{|b|\le B} \chi_2(b)\sum\limits_{i=1}^{I_{p}} 
\chi_3(u_{p,i})\chi_4(v_{p,i}) \sum\limits_{j=1}^{I_{q}}
\chi'_3(u_{q,j})\chi'_4(v_{q,j}),
$$
where $\chi_1$, $\chi_2$ are characters modulo $pq$ such that exactly one of them
is the principal character, $\chi_3$, $\chi_4$ are characters modulo $p$, and
$\chi'_3$, $\chi'_4$ are characters modulo $q$. Here the characters 
$\chi_{3,4}$, $\chi'_{3,4}$
depend on the characters $\chi_{1,2}$. Now Lemma~\ref{PolyaVinolemma} implies that
$$
E_1\ll \frac{I_{p}I_{q}(A+B)}{\sqrt{pq}}\log pq.
$$ 
 
{\it Estimation of $E_2$.}\ \ \ \  Given $k,l\in \intz$ and a character
$\chi_1$ modulo $pq$, 
the number of characters $\chi$ modulo $pq$ such that  
$(\cdot/p)_4^{k} (\cdot/q)_4^{l}(\chi\chi')^3=\chi_1$  
is $\le 9$, and the number of $\chi$ modulo $pq$ such that 
$\chi^2=\chi_1$ is $\le 4$. 
Thus, using the Cauchy-Schwarz inequality, 
we deduce that
\begin{equation} \label{cha2}
\begin{split}
E_2 &\ll \frac{1}{pq} \sum\limits_{k=1}^4 \sum\limits_{l=1}^4 \left(
\sum\limits_{\chi} \left|\sum\limits_{i=1}^{I_{p}} \left(\frac{u_{p,i}}{p}\right)_4^{k} \chi(u_{p,i}^{-3}v_{p,i}^2)
\right|^2 \right)^{\frac{1}{2}} \\
& \hspace*{.25in}\times \left(\sum\limits_{\chi'} \left|\sum\limits_{j=1}^{I_{q}} 
\left(\frac{u_{q,j}}{q}\right)_4^{l} \chi'(u_{q,j}^{-3}v_{q,j}^2)
\right|^2 \right)^{\frac{1}{2}}
\left(\sum\limits_{\chi_1\not= \chi_0\chi'_0}
\left| \sum\limits_{|a|\le A} \chi_1(a)
\right|^4\right)^{\frac{1}{4}} \left(\sum\limits_{\chi_2\not=\chi_0\chi'_0} 
\left|\sum\limits_{|b|\le B}  \chi(b)\right|^4\right)^{1/4},
\end{split}
\end{equation}  
where $\chi$ runs over all characters modulo $p$, $\chi'$ runs over all
characters modulo $q$, and $\chi_1,\chi_2$ run over all non-principal
characters modulo $pq$. \newline

By Lemma~\ref{charlemma} (i), the number of $i$'s such that 
$u_{p,i}^{-3}v_{p,i}^2$ lie in a fixed residue class modulo $p$ is bounded by 4.
The same is true for the number of $j's$ such that 
$u_{q,j}^{-3}v_{q,j}^2$ lie in a fixed residue class modulo $q$.  
Using this, Lemma~\ref{Dirlemma} and Lemma~\ref{fourthmomentlemma}, the expression on the right-hand side of
(\ref{cha2}) is dominated by 
$$
\ll (I_{p}I_{q}AB)^{1/2}\log^3 pq.
$$

{\it The final estimate.} Combining all contributions, we obtain 
\begin{equation} \label{end2}
\begin{split}
& \mathop{\sum \sum}_{\substack{|a|\le A,\ |b|\le B,\ p, q \nmid ab \\ 2\sqrt{p}\alpha\le \lambda_{E(a,b)}(p)\le 2\sqrt{p}\beta \\  2\sqrt{q}\alpha\le \lambda_{E(a,b)}(q)\le 2\sqrt{q}\beta }} 1 \\
&= 4AB\sum\limits_{2\sqrt{p}\alpha\le r \le 2\sqrt{p}\beta}
\frac{H(r^2-4p)}{2p}\sum\limits_{2\sqrt{q}\alpha\le r \le 2\sqrt{q}\beta}
\frac{H(r^2-4q)}{2q} \\ 
& \hspace*{.25in} +O\left(\frac{AB(I_{p}+I_{q})}{pq}
+ABI_{p}I_{q}\left(\frac{1}{p^2q}+\frac{1}{pq^2}\right)+
(I_{p}I_{q}AB)^{1/2}\log^3 pq+
\frac{I_{p}I_{q}(A+B)}{\sqrt{pq}}\log pq\right).
\end{split}
\end{equation}
We have proved this estimate only for distinct primes $p,q$ with $p\equiv q
\equiv 1 \pmod{4}$, but the same estimate can be proved in the cases
$pq\equiv -1 \pmod{4}$ and $p\equiv q\equiv 3 \pmod{4}$ in a similar way. 
Now we obtain Theorem~\ref{TH2} from (\ref{mu}), (\ref{start}), \eqref{goal2}, 
(\ref{app}), \eqref{app1}, 
(\ref{end2}), Theorem~\ref{maintermtheorem}, Lemma \ref{HPlemma}, $I_p\le H_p\ll p$,
$F(\alpha,\beta)\ge x^{-1/2}$ and the 
prime number theorem. \newline

\noindent{\bf Acknowledgments.}  The authors would like to thank the
referee for his/her many comments.  This paper was written when the first
and second-named authors held postdoctoral fellowships at the Department
of Mathematics and Statistics at Queen's University and the Department
of Mathematics at the University of Toronto, respectively.   The authors
wish to thank these institutions for their financial support.  More in
particular, for the generous support and encouragement given to him during and
even after his stay in the University of Toronto, the second-named author owes a debt of gratitude to
Prof. John B. Friedlander.

\vspace*{.7cm}

\noindent Department of Mathematics and Statistics, Queen's University \newline
University Ave, Kingston, ON K7L 3N6 Canada \newline
Email: {\tt sbaier@mast.queensu.ca} \newline

\noindent Department of Mathematics, University of Toronto \newline
40 Saint George Street, Toronto, ON M5S 2E4 Canada \newline
Email: {\tt lzhao@math.toronto.edu}

\end{document}